\documentclass[a4paper,12pt]{article}
\usepackage{amssymb}
{}

\newtheorem{thm}{Theorem}[section]

\newcommand{\R}{\mathbb{R}}
\renewcommand{\t}{\tilde}
\renewcommand{\r}{\rho}

\renewcommand{\O}{\Omega}
\renewcommand{\d}{\delta}
\newcommand{\p}{\partial}

\newcommand{\bc}{\begin{cor}}
\newcommand{\ec}{\end{cor}}
\newcommand{\bl}{\begin{lem}}
\newcommand{\el}{\end{lem}}
\newcommand{\bp}{\begin{prop}}
\newcommand{\ep}{\end{prop}}
\newcommand{\bt}{\begin{thm}}
\newcommand{\et}{\end{thm}}
\newcommand{\bal}{\begin{array}{ll}}
\newcommand{\ba}{\begin{array}}
\newcommand{\bac}{\begin{array}{ccc}}
\newcommand{\ea}{\end{array}}

\newcommand{\be}{\begin{equation}}
\newcommand{\ee}{\end{equation}}

\newtheorem{lem}[thm]{Lemma}
\newtheorem{cor}[thm]{Corollary}
\newtheorem{prop}[thm]{Proposition}

\begin{document}
\title {\Large\bf{Well-posedness of the IBVP for 2-D Euler Equations
with Damping
   {\thanks {The research was supported in part by NNSF of
China (10531020).}}}}

\author { Yongqin Liu{\footnote{email:~yqliu2@yahoo.com.cn}}\\
 {\small\emph{School of Mathematical Sciences, Fudan University, Shanghai, China}}\\
 Weike Wang{\footnote{email:~wkwang@sjtu.edu.cn}}\\
 \small\emph{Department of Mathematics, Shanghai Jiao Tong University, Shanghai, China}}
\date{}
\maketitle
 {\bf Abstract}. In this paper we focus on the initial-boundary value problem of the 2-D isentropic Euler
 equations with damping. We prove the global-in-time existence of classical solution
 to the initial-boundary value problem
   by
 the  method of energy estimates.\\
\indent {\bf keywords}: Euler equation; initial-boundary value problem; well-posedness. \\
\indent {\it MSC(2000)}: 35A05;~35L45.
\section{Introduction}
In this paper we concern the global-in-time well-posedness of
solutions to the initial-boundary value problem (IBVP) of the
following isentropic Euler equations with damping in two dimensional
space. \be\label{1a}
 \left\{\begin{array}{ll}
&\t{\r}_t+{\rm{div}}(\t{\r}\t{u})=0,\ \t{x}> st,\ \t{y}\in\R,\ t>0,\\&\\
&(\t{\r}\t{u}_j)_t+{\rm{div}}(\t{\r}\t{u}\t{u}_j)+P(\t{\r})_{\t{x}_j}=-k\t{\r}\t{u}_j,\
j=1,2.
\end{array}
\right. \ee   Here $\t{u}(\t{x},t)=(\t{u}_1,\t{u}_2)(\t{x},t),\
\t{\r}(\t{x},t),$ $P=P(\t{\r})$ represent the velocity, fluid
density and pressure respectively, $k>0$ is a positive constant, $s$
is a real number. As is well-known, (\ref{1a}) in one-dimension can
be written into the p-system with damping in the Lagrangian
coordinates, \be\label{1b}\left\{\begin{array}{ll}
&v_t-u_x=0,\ x\in\R^+,\ t>0,\\&\\
&u_t+P(v)_x=-ku,\ k>0.
\end{array}
\right. \ee Here $v(x,t)>0$ and $u(x,t)$ represent the specific
volume and velocity, respectively; the pressure $P(v)$ is assumed to
be a smooth function of $v$ with $P(v)>0,\ P^{\prime}(v)<0.$ In
\cite{NY} Nishihara and Yang studied the boundary effect on the
asymptotic behavior of solution to (\ref{1b}) with the Dirichlet
boundary condition $u|_{x=0}=0.$ In \cite{WY} Wang and Yang
considered the time-asymptotic behavior of solutions to the Cauchy
problem for the isentropic Euler equations with damping in
multi-dimensions, where the global existence and pointwise estimates
of the solutions are obtained, moreover they also obtained the
optimal $L^p(1<p\leq\infty)$ convergence rate of the solution when
it is a perturbation of a constant state. Moreover, in
\cite{HMS,MM,MN1,MN2}), Matsumura, et al studied the viscous shock
wave and the  asymptotic behaviors of solutions to the IBVP of the
p-system with viscosity. For the IBVP of the Navier-Stokes
equations, there are some results. In \cite{KK1} Kagei and Kobayashi
studied the large-time behavior of solutions to the compressible
Navier-Stokes equations in the half space in $\R^3.$ In \cite{KK}
Kagei and Kawashima studied the stability of planar stationary
solutions to the IBVP of the Navier-Stokes equations on the half
space. However there are few works on the IBVP in multi-dimensional
case to the Euler equations with linear damping (\ref{1a}) as far as
we know.

As for the IBVP, how to give the appropriate boundary conditions,
which is a key point to close the energy estimates, is a difficulty
we meet with, since the IBVP may be ill-posed under some boundary
conditions (see \cite{Kreiss}). What and how many boundary
conditions to give are two problems we have to solve at first.
Because the increase of the spatial dimensions and the number of the
equations, we can not simply propose the Dirichlet condition on the
velocity as in one-dimensional case (see \cite{NY}). By
 diagonalizing  the coefficient matrix of the normal (with respect to the boundary)
derivative  of the unknown variables, we give the boundary
conditions on the linear combination of the unknown variables, and
find that the number of the boundary conditions to give is
determined by the number of the positive eigenvalues of the
coefficient matrix of the normal (with respect to the boundary)
derivative  of the unknown variables.

A matter worthy of note is that the process of making a priori
estimates for IBVP is more complex than that for Cauchy problem.
Especially in dealing with the boundary terms composed of the
higher-order normal derivatives, we have to take the original system
into consideration. Moreover, the complexity increases as the order
of the derivatives grows higher. In order to close the energy
estimates, we make use of some techniques in dealing with the
boundary terms.

Another matter to mention is about the local existence of solutions.
In general, for Cauchy problem of symmetric hyperbolic systems, the
local existence of  classical solutions could be obtained  without
the assumption of small initial data (see \cite{S}), while for IBVP,
there is some difference. Since the boundary terms could affect the
symmetric structure  of the system in the process of making energy
estimates, there exists some difficulty (essential or technical) in
obtaining the local existence of solutions without the assumption of
small initial data. However, this does not affect our ultimate
results, because the global-in-time a priori estimates require that
the initial data be small. So what we need is to prove the local
existence of classical solutions in the case of small initial data,
and this could be obtained by using the iterative scheme.

The rest of the paper is as follows. After we state the notations,
in section 2 we give the a priori estimates by energy methods . In
section 3 we give the main theorems and show the global existence of
the classical solution to the IBVP.

Notations. We denote generic constants by $C$.
$\p^k\triangleq(\p_x^k,\p_x^{k-1}\p_y,\cdots,\p_y^k).$
$\O_t\triangleq\R^+\times\R\times[0,t].$ $L^p(1\leq p\leq\infty)$ is
the usual Lebesgue space with the norm $|\cdot|_p$, $W^{m,~p},m\in
{\mathbb{Z}^+},p\in[1,\infty]$ denotes the usual Sobolev space with
its norm
$$
\|f\|_{W^{m,~p}}\triangleq(\sum\limits_{k=0}\limits^{m}|\partial_x^kf|_p^p)^{1\over
p}.
$$
In particular, we use $W^{m,~2}=H^m$ with its norm $\|\cdot\|_m, and
\|\cdot\|_0=\|\cdot\|.$ Since we cope with the initial-boundary
value problem, for convenience, we denote,
$\|f\|^2(0,\cdot,t)\triangleq\int_{\R}|f(0,y,t)|^2dy,\
 \|f\|^2(t)\triangleq\int_{\R}\int_{\R^+}|f(x,y,t)|^2dxdy.$
\section{Energy estimates}
In this paper we consider the small perturbation near the constant
state $(\r^{\sharp},u^{\sharp})$, without loss of generality we
choose $\r^{\sharp}=1,u^{\sharp}=0$. The real number $s$ play an
important role in proposing the appropriate boundary conditions. The
comparison between $s$ and $r$ decides the number of the boundary
conditions we could propose. In this paper we consider the case
$0<s<r$, and the other cases can be studied in the future.
Correspondingly we study the following initial-boundary value
problem,

\be \label{3a} \left\{ \bal&\t{\r}_t+{\rm{div}}(\t{\r}\t{u})=0,\
\t{x}> st,\ \t{y}\in\R,\
t>0,\\&\\&(\t{\r}\t{u}_1)_t+{\rm{div}}(\t{\r}\t{u}\t{u}_1)
+P(\t{\r})_{\t{x}}=-k\t{\r}\t{u}_1,\\&\\&(\t{\r}\t{u}_2)_t+{\rm{div}}(\t{\r}\t{u}\t{u}_2)
+P(\t{\r})_{\t{y}}=-k\t{\r}\t{u}_2,
\\&\\&(\t{\r},\t{u}_1,\t{u}_2)(\t{x},\t{y},t)|_{t=0}=(\r_0+1,ru_{10},ru_{20})(\t{x},\t{y}),\\&\\&
(\t{\r}+{{\t{u}_1}\over r})|_{\t{x}=st}=1,\ea\right.
 \ee
 where
$r^2=P^{\prime}(1)>0$, $\r_0,u_{10},u_{20}$ are given functions, and
$\inf\limits_{(x,y)\in \R^+\times\R}\r_0(x,y)+1>0.$ We assume that
the pressure $P(\t\r)$ is smooth in a neighborhood of
$\r^{\sharp}=1$.

Next we will make a series of transformations to the coordinates and
unknown variables. First, $\t{x}\rightarrow x+st,\ t\rightarrow t,$
changes the domain we study  from a wedge to the half space. Second,
the translation transformation $\t{\r}\rightarrow \bar{\r}-1,
\t{u}_1\rightarrow \bar{u}_1,\ \t{u}_2\rightarrow \bar{u}_2,$
  linearizes (\ref{3a}). Last, the scaling
 transformation $\bar{\r}\rightarrow {\r}, \bar{u}_1\rightarrow
r{u_1},\ \bar{u}_2\rightarrow r{u_2},$  reformulates the problem
(\ref{3a}) to the following system,  \be \label{3b} \left\{ \bal&
\r_t-s\r_x+ru_{1x}+ru_{2y}=-r{\rm{div}}(\r u) ,\ x>0,\ y\in\R,\ t>0,
\\&\\&u_{1t}-su_{1x}+r\r_{x}+ku_1=-ru\cdot\nabla u_1+{1\over r}(r^2-{{P^{\prime}(1+\r)}\over{1+\r}})\r_x,
\\&\\&u_{2t}-su_{2x}+r\r_{y}+ku_2=-ru\cdot\nabla u_2+{1\over r}(r^2-{{P^{\prime}(1+\r)}\over{1+\r}})\r_y,
\\&\\&(\r,u_1,u_2)(x,y,t)|_{t=0}=(\r_0,u_{10},u_{20})(x,y),
\\&\\&(\r+u_1)|_{x=0}=0,
 \ea\right.
 \ee where $u=(u_1,u_2)$.
Denote $B=r^2-{{P^{\prime}(1+\r)}\over{1+\r}},$ $$ A_1=\left( \bac
-s&r&0\\&&\\ r&-s&0\\&&\\0& 0&-s \ea \right), A_2=\left( \bac
0&0&r\\&&\\ 0&0&0\\&&\\r& 0&0 \ea \right), A_3=\left( \bac
0&0&0\\&&\\ 0&k&0\\&&\\0& 0&k \ea \right),
$$
$$
H=\left(\bal&h_1\\&\\&h_2\\&\\&h_3\ea\right)=\left(\bal&-r{\rm{div}}(\r
u)\\&\\&-ru\cdot\nabla u_1+{B\over r}\r_x\\&\\&-ru\cdot\nabla
u_2+{B\over r}\r_y\ea\right),\
W=\left(\bal&\r\\&\\&u_1\\&\\&u_2\ea\right).
$$
Then we can rewrite (\ref{3b}) as following,
$$W_t+A_1W_x+A_2W_y+A_3W=H.$$ In order to diagonalize the coefficient matrix $A_1$, we introduce an orthogonal
transform.
Let $$ S_0=\left( \bac {\sqrt{2}\over2}&{\sqrt{2}\over2}&0\\&&\\
{\sqrt{2}\over2}&-{\sqrt{2}\over2}&0\\&&\\0& 0&1 \ea \right),
$$ then $S_0^{-1}=S_0.$
Denote $S_i=S_0A_iS_0, \ i=1,2,3, \ V=S_0W,$ then
$$
S_1=\left( \bac -s+r&0&0\\&&\\ 0&-s-r&0\\&&\\0& 0&-s \ea
\right),\  S_2=\left( \bac 0&0&{\sqrt{2}\over2}r\\&&\\
0&0&{\sqrt{2}\over2}r\\&&\\{\sqrt{2}\over2}r& {\sqrt{2}\over2}r&0
\ea \right),
$$$$
\ S_3=\left( \bac {k\over2}&-{k\over2}&0\\&&\\
-{k\over2}&{k\over2}&0\\&&\\0& 0&k \ea \right),\
V=\left(\ba{cc}&{\sqrt{2}\over2}(\r+u_1)\\&\\&{\sqrt{2}\over2}(\r-u_1)\\&\\&u_2\ea\right)
\triangleq\left(\bal&v_1\\&\\&v_2\\&\\&v_3\ea\right),
$$
$V(x,y,0)\triangleq V_0(x,y),$ thus we reformulate (\ref{3b}) to the
following problem, \be\label{3c} \left\{ \bal&
V_t+S_1V_x+S_2V_y+S_3V=S_0H,\\&\\&
V(x,y,0)=V_0(x,y),\\&\\&v_1(x,y,t)|_{x=0}=0.
 \ea
\right. \ee

Specifically, (\ref{3c}) can be written into the following form,
\be\label{3d} \left\{ \bal&
v_{1t}-(s-r)v_{1x}+{\sqrt{2}\over2}rv_{3y}+{k\over2}(v_1-v_2)={\sqrt{2}\over2}(h_1+h_2),
 \\&\\&v_{2t}-(s+r)v_{2x}+{\sqrt{2}\over2}rv_{3y}+{k\over2}(v_2-v_1)={\sqrt{2}\over2}(h_1-h_2),\\&\\
 &v_{3t}-sv_{3x}+{\sqrt{2}\over2}r(v_{1y}+v_{2y})+kv_3=h_3,\\&\\&
(v_1,v_2,v_3)(x,y,0)=({\sqrt{2}\over2}(\r_0+u_{10}),{\sqrt{2}\over2}(\r_0-u_{10}),u_{20})(x,y),\\&\\&v_1(x,y,t)|_{x=0}=0.
 \ea
\right. \ee

In the following we will estimate $(\r,u_1,u_2)$ under the a priori
assumption \be\label{a}
N(T)\triangleq\sup\limits_{0<t<T}\{\|W\|^2_l(t)\}\leq \d_0, \
0<\d_0\ll1,~l\geq4. \ee  By Sobolev inequality and the system
(\ref{3b}), we know that
$$
\sum\limits_{0\leq k_1+k_2+k_3\leq
l-2}\sup\limits_{\O_T}|\p^{k_1}_x\p^{k_2}_{y}\p_t^{k_3}W|\leq
C\d_0,$$$$|B|=|r^2-{{P^{\prime}{(1+\r)}}\over{1+\r}}|\leq C|\r|\leq
C\d_0.
$$
Now we will obtain a series of estimates corresponding to the
k-order derivatives  (k=0,1,2,3,4), denoted by Estimate A, B, C, D
and E, and higher order derivatives of the solution in order to
close the energy estimates . In the process of energy estimates, we
use the fact that $\|\p^kV\|=\|\p^kW\|$, $k\geq0$ is an integer,
since $S_0$ is an orthogonal matrix.

\subsection{Estimate A}
Multiplying (\ref{3c}) by $V$ and integrating it over $\O_t$, since
$$
|\int^t_0\int^{\infty}_{-\infty}\int^{\infty}_0H\cdot
Wdxdyd\tau|\leq
C\d_0\int^t_0(\|\r\|^2(0,\cdot,\tau)+\|u\|^2_1(\tau)+\|\nabla\r\|^2(\tau))d\tau,
$$
we get that \be\label{2101} \bal&\ \ \
\|W\|^2(t)+\int^t_0(\|W\|^2(0,\cdot,\tau))+\|u\|^2(\tau))d\tau\\&\\&\leq
C\|W_0\|^2+C\d_0\int^t_0(\|u\|^2_1(\tau)+\|\nabla\r\|^2(\tau))d\tau.\ea
\ee
\subsection{Estimate B}
By direct calculation we obtain the estimates on the nonlinear
terms. \bl\label{22a} Assume (\ref{a}) holds, then
$$
\bal&\ \ \ |\int^t_0\int^{\infty}_{-\infty}\int^{\infty}_0\p_yH\cdot
\p_yWdxdyd\tau|\\&\\&\leq
C\d_0(\|\p_y\r\|^2(t)+\|\p_y\r_0\|^2\\&\\&\ \ \
+\int^t_0[\|W\|^2_1(0,\cdot,\tau)+\|u\|_1^2(\tau)+\|\nabla\r\|^2(\tau)]d\tau),\ea
$$

$$
\bal&\ \ \ |\int^t_0\int^{\infty}_{-\infty}\int^{\infty}_0\p_xH\cdot
\p_xWdxdyd\tau|\\&\\&\leq
C\d_0(\|\p_x\r\|^2(t)+\|\p_x\r_0\|^2\\&\\&\ \ \
+\int^t_0[\|W\|^2_1(0,\cdot,\tau)+\|u\|_1^2(\tau)+\|\nabla\r\|^2(\tau)]d\tau).\ea
$$
 \el

As for the boundary terms we have the following estimates.
\bl\label{22b} Assume (\ref{a}) holds, then \be\label{2201}
\|\p_xv_1\|(0,\cdot,t)\leq
C(\|W\|^2+\|\p_yu_2\|^2)(0,\cdot,t)+C\d_0\|W\|_1^2(0,\cdot,t).\ee
\el

{\bf Proof.} By virtue of $(\ref{3d})_1$, we get that
$$\bal
\|\p_xv_{1}\|^2(0,\cdot,t)&\leq
C(\|\p_yv_{3}\|^2+\|v_2\|^2+\|h_1+h_2\|^2)(0,\cdot,t)\\&\\&\leq
(\|W\|^2+\|\p_yu_2\|^2)(0,\cdot,t)+C\d_0\|W\|_1^2(0,\cdot,t).\ea
$$ Thus (\ref{2201}) is proved. $\Box$

  Multiplying $\p_y(\ref{3c})$ by $\p_yV$ and integrating it over $\O_t$, combined with lemma \ref{22a}, yields
that, \be\label{2202}\bal&\ \ \
\|\p_yW\|^2(t)+\int^t_0(\|\p_yW\|^2(0,\cdot,\tau)+\|\p_yu\|^2(\tau))d\tau\\&\\&\leq
C\|\p_yW_0\|^2+C\d_0\int^t_0[\|\p_yW\|^2(0,\cdot,\tau)+(\|\nabla\r\|^2+\|
u\|^2_1)(\tau)]d\tau. \ea\ee

  Multiplying $\p_x(\ref{3c})$ by $\p_xV$ and integrating it over $\O_t$, combined with lemma \ref{22a}, yields
that, \be\label{2203}\bal&\ \ \
\|\p_xW\|^2(t)+\int^t_0[(\|\p_x(\r-u_{1})\|^2+\|\p_xu_{2}\|^2)(0,\cdot,\tau)+
\|\p_xu\|^2(\tau)]d\tau\\&\\&\leq
C(\|\p_xW_0\|^2+\int^t_0\|\p_xv_1\|^2(0,\cdot,\tau)d\tau)\\&\\&\ \ \
+C\d_0\int^t_0[\|\p_xW\|^2(0,\cdot,\tau)+\|\nabla\r\|^2(\tau)+\|
u\|^2_1(\tau)]d\tau. \ea \ee

 Choose $\lambda_1$ suitably small such
that $(\ref{2202})+\lambda_1 (\ref{2203})$, combined with
(\ref{2201}) and (\ref{2101})
 yields that
\be\label{2205}\bal&\|\p
W\|^2(t)+\int^t_0[(\|\p_yW\|+\|\p_x(\r-u_1)\|^2+\|\p_xu_2\|^2)(0,\cdot,\tau)
+\|\p u\|^2(\tau)]d\tau\\&\\&\leq C(\|W_0\|_1^2+C\d_0\int^t_0[\|\p
W\|^2(0,\cdot,\tau)+(\|\nabla\r\|+\|u\|_1^2)(\tau)]d\tau.\ea \ee

 (\ref{2205}),~(\ref{2201})
and (\ref{2101}) yield that
 \be\label{2206}
\|W\|^2_1(t)+\int^t_0[\|W\|^2_1(0,\cdot,\tau)+\|u\|^2_1(\tau)]d\tau\leq
C\|W_0\|^2_1+C\d_0\int^t_0\|\nabla\r\|^2(\tau)d\tau. \ee

From (\ref{3b}), we have that $$\|W_t\|^2(0,\cdot,t)\leq
C\|W\|^2_1(0,\cdot,t),\ \  \|W_t\|^2(t)\leq C\|W\|^2_1(t).
$$

Thus (\ref{2206}) yields that \be \label{2207}\bal &\ \ \
\|W\|^2_1(t)+
\|W_t\|^2(t)+\int^t_0[(\|W\|^2_1+\|W_t\|^2)(0,\cdot,\tau)+\|u\|^2_1(\tau)]d\tau\\&\\&\leq
C\|W_0\|^2_1+C\d_0\int^t_0\|\nabla\r\|^2(\tau)d\tau. \ea\ee

Since
$$u_{1t}\r_x=(u_1\r_x)_t-(u_1\r_t)_x+u_{1x}\r_t,\ \
u_{2t}\r_y=(u_2\r_y)_t-(u_2\r_t)_y+u_{2y}\r_t,$$ by virtue of Cauchy
inequality, $(\ref{3b})_1\r_t+s(\ref{3b})_2\r_x+(\ref{3b})_3\r_y$
yields that
\be\label{2208}\bal&\int^t_0(\|\r_t\|^2(\tau)+\|\nabla\r\|^2(\tau))d\tau\leq\\&\\&
C(\|W_0\|^2_1+(\|\nabla\r\|^2+\|u\|^2)(t)+
\int^t_0[(\|\r_t\|^2+\|u\|^2)(0,\cdot,\tau)+\|u\|_1^2(\tau)]d\tau)\\&\\&+C\d_0\int_0^t(\|u\|_1^2
+\|\nabla\r\|^2)(\tau)d\tau. \ea\ee

Choose $\lambda_2$ suitably small such that
$(\ref{2207})+\lambda_2(\ref{2208})$ yields that
 \be\label{2204}\bal&
\|W\|^2_1(t)+\|W_t\|^2(t)
+\int^t_0[(\|W\|^2_1+\|W_t\|^2)(0,\cdot,\tau)\\&\\&
+(\|\r_t\|^2+\|\nabla\r\|^2+\|u\|^2_1) (\tau)]d\tau\leq
C\|W_0\|^2_1. \ea\ee

From (\ref{3b}) we know that $\|u_t\|^2(t)\leq
C(\|\nabla\r\|^2+\|u\|^2_1)(t)$, thus (\ref{2204}) yields that
 \be\label{2209}\bal&
\|W\|^2_1(t)+\|W_t\|^2(t)
+\int^t_0[(\|W\|^2_1+\|W_t\|^2)(0,\cdot,\tau)\\&\\&+(\|W_t\|^2+\|\nabla\r\|^2+\|u\|^2_1)
(\tau)]d\tau\leq C\|W_0\|^2_1. \ea\ee
\subsection{Estimate C}
By direct calculation we have the following estimates on the
nonlinear terms.
\bl\label{23a} Assume (\ref{a}) holds, then
$$\bal&\ \
 |\int^t_0\int^{\infty}_{-\infty}\int^{\infty}_0\p_y^2H\cdot
\p_y^2Wdxdyd\tau|\\&\\&\leq
C\d_0(\|\p_y^2\r_0\|^2+\|\p_y^2\r\|^2(t)\\&\\&\ \
+\int^t_0[\|\p^2W\|^2(0,\cdot,\tau)+(\|\nabla\r\|^2_1+\|u\|^2_2)(\tau)]d\tau),
\ea
$$

$$\bal&\ \
 |\int^t_0\int^{\infty}_{-\infty}\int^{\infty}_0\p_y\p_xH\cdot
\p_y\p_xWdxdyd\tau|\\&\\&\leq
C\d_0(\|\p_y\p_x\r_0\|^2+\|\p_y\p_x\r\|^2(t)\\&\\&\ \
+\int^t_0[\|\p^2W\|^2(0,\cdot,\tau)+(\|\nabla\r\|^2_1+\|u\|^2_2)(\tau)]d\tau),
\ea
$$

$$\bal&\ \
 |\int^t_0\int^{\infty}_{-\infty}\int^{\infty}_0\p_x^2H\cdot
\p_x^2Wdxdyd\tau|\\&\\&\leq
C\d_0(\|\p_x^2\r_0\|^2+\|\p_x^2\r\|^2(t)\\&\\&\ \
+\int^t_0[\|\p^2W\|^2(0,\cdot,\tau)+(\|\nabla\r\|^2_1+\|u\|^2_2)(\tau)]d\tau).
\ea
$$
 \el

 As for the boundary terms we have the following estimates.
 \bl\label{23b} Assume (\ref{a}) holds, then
 \be\label{2305}
\|\p_y\p_xv_1\|^2(0,\cdot,t)\leq
C(\|W\|_1^2+\|\p_y^2u_2\|)(0,\cdot,t)+C\d_0\|W\|^2_2(0,\cdot,t),\ee
\be\label{2308} \|\p_x^2v_1\|^2(0,\cdot,t)\leq
C(\|W\|_1^2+\|\p_y^2(\r-u_1)\|+\|\p_y\p_xu_2\|^2)(0,\cdot,t)+C\d_0\|W\|^2_2(0,\cdot,t).\ee
 \el
{\bf Proof.} In view of $\p_y(\ref{3d})_1$, we get that
$$
\bal\|\p_y\p_xv_1\|^2(0,\cdot,t)&\leq
C(\|\p_y^2v_3\|^2+\|\p_y(v_1-v_2)\|^2+\|\p_y(h_1+h_2)\|^2)(0,\cdot,t)\\&\\&\leq
C(\|W\|_1^2+\|\p_y^2u_2\|)(0,\cdot,t)+C\d_0\|W\|^2_2(0,\cdot,t).\ea
$$
Thus (\ref{2305}) is proved.

 In view of $\p_x(\ref{3d})_1,$ we get
that
$$\bal&\ \ \ \|\p_x^2v_1\|^2(0,\cdot,t)\\&\\&\leq
C(\|\p_x\p_tv_1\|+\|\p_y\p_xv_3\|^2+\|\p_x(v_1-v_2)\|^2+\|\p_x(h_1+h_2)\|^2)(0,\cdot,t)\\&\\&\leq
C\|\p_x\p_tv_1\|(0,\cdot,t)+C(\|W\|_1^2+\|\p_y\p_xu_2\|^2)(0,\cdot,t)+C\d_0\|W\|^2_2(0,\cdot,t)
\\&\\&\leq C\|\p_y\p_tv_3\|(0,\cdot,t)+C(\|W\|_1^2+\|\p_y\p_xu_2\|^2)(0,\cdot,t)+C\d_0\|W\|^2_2(0,\cdot,t)
\\&\\&\leq
C(\|W\|_1^2+\|\p_y^2(\r-u_1)\|+\|\p_y\p_xu_2\|^2)(0,\cdot,t)+C\d_0\|W\|^2_2(0,\cdot,t).\ea
$$
Thus (\ref{2308}) is proved. \ \ $\Box$

Multiplying $\p_y^2(\ref{3c})$ by $\p_y^2V$ and integrating it over
$\O_t$, we have that \be\label{2301}\bal&\ \
\|\p_y^2W\|^2(t)+\int^t_0[\|\p_y^2W\|^2(0,\cdot,\tau)+\|\p^2_yu\|^2(\tau)]d\tau\\&\\&\leq
C(\|\p_y^2W_0\|^2+|\int^t_0\int^{\infty}_{-\infty}\int^{\infty}_0\p_y^2H\cdot
\p_y^2Wdxdyd\tau|). \ea\ee

Similarly, we have that
 \be\label{2303}\bal&\ \ \
\|\p_y\p_xW\|^2(t)+\int^t_0[(\|\p_y\p_x(\r-u_1)\|^2+\|\p_y\p_xu_2\|^2)(0,\cdot,\tau)\\&\\&\
\ \ +\|\p_y\p_xu\|^2(\tau)] d\tau\\&\\&\leq
C(\|\p_y\p_xW_0\|^2+\int^t_0\|\p_y\p_xv_1\|^2(0,\cdot,\tau)d\tau\\&\\&\
\  \  +|\int^t_0\int^{\infty}_{-\infty}\int^{\infty}_0
\p_y\p_xH\cdot \p_y\p_xWdxdyd\tau|), \ea\ee

and
 \be\label{2307}\bal&\ \ \
\|\p_x^2W\|^2(t)+\int^t_0[(\|\p_x^2(\r-u_1)\|^2+\|\p_x^2u_2\|^2)
(0,\cdot,\tau)+\|\p_x^2u\|^2(\tau)] d\tau\\&\\&\leq
C(\|\p_x^2W_0\|^2+\int^t_0\|\p_x^2v_1\|^2(0,\cdot,\tau)d\tau+|\int^t_0\int^{\infty}_{-\infty}\int^{\infty}_0
\p_x^2H\cdot \p_x^2Wdxdyd\tau|). \ea\ee

Choose $\lambda_3,~\lambda_4$ suitably small such that
$(\ref{2301})+\lambda_3(\ref{2303})$, combined with (\ref{2305}) and
(\ref{2209}), yields that \be\label{2306}\bal&\ \
\|\p_y^2W\|^2(t)+\|\p_y\p_xW\|^2(t)+\int^t_0[(\|\p_y^2W\|^2+\|\p_y\p_x(\r-u_1)\|^2\\&\\&\
\  +\|\p_y\p_xu_2\|^2)(0,\cdot,\tau)
+(\|\p_y^2u\|^2+\|\p_y\p_xu\|^2)(\tau)] d\tau\\&\\&\leq
C(\|W_0\|_2^2+\d_0\int^t_0\|W\|_2^2(0,\cdot,\tau)d\tau\\&\\&\ \
+|\int^t_0\int^{\infty}_{-\infty}\int^{\infty}_0 \p_y^2H\cdot
\p_y^2Wdxdyd\tau|+|\int^t_0\int^{\infty}_{-\infty}\int^{\infty}_0\p_y\p_xH\cdot
\p_y\p_xWdxdyd\tau|), \ea\ee

and $(\ref{2306})+\lambda_4(\ref{2307})$, combined with (\ref{2308})
and (\ref{2209}), yields that \be\label{2309}\bal&\ \
\|\p^2W\|^2(t)+\int^t_0[(\|\p_y^2W\|^2+\|\p_y\p_x(\r-u_1)\|^2\\&\\&\
\ +
\|\p_y\p_xu_2\|^2+\|\p_x^2(\r-u_1)\|^2+\|\p_x^2u_2\|^2)(0,\cdot,\tau)
+\|\p^2u\|^2(\tau)] d\tau\\&\\&\leq
C(\|W_0\|_2^2+\d_0\int^t_0\|W\|_2^2(0,\cdot,\tau)d\tau+|\int^t_0\int^{\infty}_{-\infty}\int^{\infty}_0\p_y^2H\cdot
\p_y^2Wdxdyd\tau|\\&\\&\ \
+|\int^t_0\int^{\infty}_{-\infty}\int^{\infty}_0 \p_y\p_xH\cdot
\p_y\p_xWdxdyd\tau| +|\int^t_0\int^{\infty}_{-\infty}\int^{\infty}_0
\p_x^2H\cdot \p_x^2Wdxdyd\tau|). \ea\ee

Combined with (\ref{2209}), lemma \ref{23a} and lemma \ref{23b},
(\ref{2309}) yields that
 \be\label{2312}\bal&\ \ \
\|W\|_2^2(t)+\int^t_0[\|W\|^2_2(0,\cdot,\tau)+\|u\|^2_2(\tau)]d\tau\\&\\&\leq
C\|W_0\|^2_2+C\d_0\int^t_0\|\nabla\r\|_1(\tau)d\tau.\ea \ee

From (\ref{3b}) it is easy to know that
$$
\|W_t\|_1^2(0,\cdot,t)\leq C\|W\|^2_2(0,\cdot,t),\ \
\|W_t\|_1^2(t)\leq C\|W\|_2^2(t),
$$
so (\ref{2312}) yields that
\be\label{2313}\bal&\ \  \
\|W\|_2^2(t)+\|W_t\|_1^2(t)+\int^t_0[(\|W\|^2_2+\|W_t\|_1^2)(0,\cdot,\tau)+\|u\|^2_2(\tau)]d\tau\\
&\\&\leq C\|W_0\|^2_2+C\d_0\int^t_0\|\nabla\r\|_1(\tau)d\tau. \ea\ee

By similar calculation to (\ref{2208}),
$\nabla(\ref{3b})_1\nabla\r_t+s\nabla(\ref{3b})_2\nabla\r_x+\nabla(\ref{3b})_3\nabla\r_y$
yields that
\be\label{2314}\bal&\int^t_0(\|\nabla\r_t\|^2+\|\nabla\r_y\|^2+\|\nabla\r_x\|^2)(\tau)d\tau\leq\\&\\&
C(\|W_0\|^2_2+(\|\nabla\r\|_1^2+\|u\|_1^2)(t)+
\int^t_0[(\|\r_t\|_1^2+\|W\|_1^2)(0,\cdot,\tau)+\|u\|_2^2(\tau)]d\tau)\\&\\&+C\d_0\int_0^t(\|u\|_2^2
+\|\nabla\r\|_1^2)(\tau)d\tau. \ea\ee

Choose $\lambda_5$ suitably small such that
$(\ref{2313})+\lambda_5(\ref{2314})$ yields that
 \be\label{2315}\bal&
\|W\|^2_2(t)+\|W_t\|_1^2(t)
+\int^t_0[(\|W\|^2_2+\|W_t\|_1^2)(0,\cdot,\tau)\\&\\&+(\|\r_t\|_1^2+\|\nabla\r\|_1^2+\|u\|^2_2)
(\tau)]d\tau\leq C\|W_0\|^2_2. \ea\ee

From (\ref{3b}) we know that $\|u_t\|_1^2(t)\leq
C(\|\nabla\r\|_1^2+\|u\|^2_2)(t)$, thus (\ref{2315}) yields that
 \be\label{2316}\bal&
\|W\|^2_2(t)+\|W_t\|_1^2(t)
+\int^t_0[(\|W\|^2_2+\|W_t\|_1^2)(0,\cdot,\tau)\\&\\&+(\|W_t\|_1^2+\|\nabla\r\|_1^2+\|u\|^2_2)
(\tau)]d\tau\leq C\|W_0\|^2_2. \ea\ee
\subsection{Estimate D}
By direct and a little tedious calculation, we get the following
estimates on the nonlinear terms.
\bl\label{24a} Assume (\ref{a})
holds, then
$$\bal&\ \
 |\int^t_0\int^{\infty}_{-\infty}\int^{\infty}_0\p_y^3H\cdot
\p_y^3Wdxdyd\tau|\\&\\&\leq
C\d_0(\|\p_y^3\r_0\|^2+\|\p_y^3\r\|^2(t)\\&\\&\ \
+\int^t_0[\|\p^3W\|^2(0,\cdot,\tau)+(\|\nabla\r\|^2_2+\|u\|^2_3)(\tau)]d\tau),
\ea
$$

$$\bal&\ \
 |\int^t_0\int^{\infty}_{-\infty}\int^{\infty}_0\p_y^2\p_xH\cdot
\p_y^2\p_xWdxdyd\tau|\\&\\&\leq
C\d_0(\|\p_y^2\p_x\r_0\|^2+\|\p_y^2\p_x\r\|^2(t)\\&\\&\ \
+\int^t_0[\|\p^3W\|^2(0,\cdot,\tau)+(\|\nabla\r\|^2_2+\|u\|^2_3)(\tau)]d\tau),
\ea
$$

$$\bal&\ \
 |\int^t_0\int^{\infty}_{-\infty}\int^{\infty}_0\p_y\p_x^2H\cdot
\p_y\p_x^2Wdxdyd\tau|\\&\\&\leq
C\d_0(\|\p_y\p_x^2\r_0\|^2+\|\p_y\p_x^2\r\|^2(t)\\&\\&\ \
+\int^t_0[\|\p^3W\|^2(0,\cdot,\tau)+(\|\nabla\r\|^2_2+\|u\|^2_3)(\tau)]d\tau),
\ea
$$

$$\bal&\ \
 |\int^t_0\int^{\infty}_{-\infty}\int^{\infty}_0\p_x^3H\cdot
\p_x^3Wdxdyd\tau|\\&\\&\leq
C\d_0(\|\p_x^3\r_0\|^2+\|\p_x^3\r\|^2(t)\\&\\&\ \
+\int^t_0[\|\p^3W\|^2(0,\cdot,\tau)+(\|\nabla\r\|^2_2+\|u\|^2_3)(\tau)]d\tau).
\ea
$$
\el

As for the estimates on the boundary terms, we have the following
results.
\bl\label{24b} Assume (\ref{a}) holds, then
\be\label{2406}\|\p_y^2\p_xv_1\|^2(0,\cdot,t)\leq
C(\|W\|_2^2+\|\p_y^3u_2\|^2)(0,\cdot,t)+C\d_0\|W\|^2_3(0,\cdot,t),\ee
\be\label{2407}\bal \|\p_y\p_x^2v_1\|^2(0,\cdot,t)&\leq
C(\|W\|_2^2+\|\p_y^2\p_xu_2\|^2+\|\p_y^3(\r-u_1)\|^2)(0,\cdot,t)\\&\\&\
\ \ +C\d_0\|W\|^2_3(0,\cdot,t),\ea\ee \be\label{2408}\bal
\|\p_x^3v_1\|^2(0,\cdot,t)&\leq
C(\|W\|_2^2+\|\p_y^3u_2\|^2+\|\p_y\p_x^2u_2\|^2\\&\\
&\ \ \ +\|\p_y^2\p_x(\r-u_1)\|^2)(0,\cdot,t)
+C\d_0\|W\|^2_3(0,\cdot,t).\ea\ee
 \el

The proof of lemma \ref{24b} is similar to that of lemma \ref{23b},
so we omit here.

Multiplying $\p_y^3(\ref{3c})$ by $\p_y^3V$ and integrating it over
$\O_t$, we have that \be\label{2409}\bal&\ \
\|\p_y^3W\|^2(t)+\int^t_0[\|\p_y^3W\|^2(0,\cdot,\tau)+\|\p^3_yu\|^2(\tau)]d\tau\\&\\&\leq
C(\|\p_y^3W_0\|^2+|\int^t_0\int^{\infty}_{-\infty}\int^{\infty}_0\p_y^3H\cdot
\p_y^3Wdxdyd\tau|). \ea\ee

Similarly, we have that
 \be\label{2410}\bal&\ \ \
\|\p_y^2\p_xW\|^2(t)\\&\\&\ \ \ +
\int^t_0[(\|\p_y^2\p_x(\r-u_1)\|^2+\|\p_y^2\p_xu_2\|^2)(0,\cdot,\tau)+\|\p_y^2\p_xu\|^2(\tau)]
d\tau\\&\\&\leq
C(\|\p_y^2\p_xW_0\|^2+\int^t_0\|\p_y^2\p_xv_1\|^2(0,\cdot,\tau)d\tau\\&\\&\
\ \  +|\int^t_0\int^{\infty}_{-\infty}\int^{\infty}_0
\p_y^2\p_xH\cdot \p_y^2\p_xWdxdyd\tau|), \ea\ee

\be\label{2411}\bal&\ \ \ \|\p_y\p_x^2W\|^2(t)\\&\\&\ \ \ +
\int^t_0[(\|\p_y\p_x^2(\r-u_1)\|^2+\|\p_y\p_x^2u_2\|^2)(0,\cdot,\tau)+\|\p_y\p_x^2u\|^2(\tau)]
d\tau\\&\\&\leq
C(\|\p_y\p_x^2W_0\|^2+\int^t_0\|\p_y\p_x^2v_1\|^2(0,\cdot,\tau)d\tau\\&\\&\
\ \  +|\int^t_0\int^{\infty}_{-\infty}\int^{\infty}_0
\p_y\p_x^2H\cdot \p_y\p_x^2Wdxdyd\tau|), \ea\ee

and
 \be\label{2412}\bal&\ \ \
\|\p_x^3W\|^2(t)+\int^t_0[(\|\p_x^3(\r-u_1)\|^2+\|\p_x^3u_2\|^2)
(0,\cdot,\tau)+\|\p_x^3u\|^2(\tau)] d\tau\\&\\&\leq
C(\|\p_x^3W_0\|^2+\int^t_0\|\p_x^3v_1\|^2(0,\cdot,\tau)d\tau+|\int^t_0\int^{\infty}_{-\infty}\int^{\infty}_0
\p_x^3H\cdot \p_x^3Wdxdyd\tau|). \ea\ee

Choose $\lambda_6,~\lambda_7,~\lambda_8$ suitably small such that
$(\ref{2409})+\lambda_6(\ref{2410})$, combined with (\ref{2406}) and
(\ref{2316}), yields that
 \be\label{2413}\bal&\ \
\|\p_y^3W\|^2(t)+\|\p_y^2\p_xW\|^2(t)+\int^t_0[(\|\p_y^3W\|^2+\|\p_y^2\p_x(\r-u_1)\|^2\\&\\&\
 \   +\|\p_y^2\p_xu_2\|^2)(0,\cdot,\tau)
+(\|\p_y^3u\|^2+\|\p_y^2\p_xu\|^2)(\tau)] d\tau\\&\\&\leq
C(\|W_0\|_3^2+\d_0\int^t_0\|W\|_3^2(0,\cdot,\tau)d\tau\\&\\&\ \
+|\int^t_0\int^{\infty}_{-\infty}\int^{\infty}_0 \p_y^3H\cdot
\p_y^3Wdxdyd\tau|+|\int^t_0\int^{\infty}_{-\infty}\int^{\infty}_0\p_y^2\p_xH\cdot
\p_y^2\p_xWdxdyd\tau|), \ea\ee

 $(\ref{2413})+\lambda_7(\ref{2411})$, combined with (\ref{2407})
and (\ref{2316}), yields that
 \be\label{2414}\bal&\ \ \
(\|\p_y^3W\|^2+\|\p_y^2\p_xW\|^2+\|\p_y\p_x^2W\|^2)(t)\\&\\&\ \ \ +
\int^t_0[(\|\p_y^3W\|^2+\|\p_y^2\p_x(\r-u_1)\|^2+
\|\p_y^2\p_xu_2\|^2+ \|\p_y\p_x^2(\r-u_1)\|^2\\&\\&\ \ \ +
\|\p_y\p_x^2u_2\|^2)(0,\cdot,\tau)
+(\|\p_y^3u\|^2+\|\p_y^2\p_xu\|^2+\|\p_y\p_x^2u\|^2)(\tau)]
d\tau\\&\\&\leq
C(\|W_0\|_3^2+\d_0\int^t_0\|W\|_3^2(0,\cdot,\tau)d\tau\\&\\&\ \ +
|\int^t_0\int^{\infty}_{-\infty}\int^{\infty}_0\p_y^3H\cdot
\p_y^3Wdxdyd\tau|+|\int^t_0\int^{\infty}_{-\infty}\int^{\infty}_0
\p_y^2\p_xH\cdot \p_y^2\p_xWdxdyd\tau|\\&\\&\ \
+|\int^t_0\int^{\infty}_{-\infty}\int^{\infty}_0 \p_y\p_x^2H\cdot
\p_y\p_x^2Wdxdyd\tau|), \ea\ee

and $(\ref{2414})+\lambda_8(\ref{2412})$, combined with (\ref{2408})
and (\ref{2316}), yields that
 \be\label{2415}\bal&\ \
\|\p^3W\|^2(t)+\int^t_0[(\|\p_y^3W\|^2+\|\p_y^2\p_x(\r-u_1)\|^2+
\|\p_y^2\p_xu_2\|^2+ \|\p_y\p_x^2u_2\|^2\\&\\&\ \ \
+\|\p_y\p_x^2(\r-u_1)\|^2+\|\p_x^3(\r-u_1)\|^2+
\|\p_x^3u_2\|^2)(0,\cdot,\tau) +\|\p^3u\|^2(\tau)] d\tau\\&\\&\leq
C(\|W_0\|_3^2+\d_0\int^t_0\|W\|_3^2(0,\cdot,\tau)d\tau\\&\\&\ \ \
+|\int^t_0\int^{\infty}_{-\infty}\int^{\infty}_0\p_y^3H\cdot
\p_y^3Wdxdyd\tau| + |\int^t_0\int^{\infty}_{-\infty}\int^{\infty}_0
\p_y^2\p_xH\cdot \p_y^2\p_xWdxdyd\tau|\\&\\&\ \ \
+|\int^t_0\int^{\infty}_{-\infty}\int^{\infty}_0 \p_y\p_x^2H\cdot
\p_y\p_x^2Wdxdyd\tau|+|\int^t_0\int^{\infty}_{-\infty}\int^{\infty}_0
\p_x^3H\cdot \p_x^3Wdxdyd\tau|). \ea\ee

Combined with lemma \ref{24a}, lemma \ref{24b} and (\ref{2316}),
(\ref{2415}) yields that \be\label{2417}\bal&\ \ \
\|W\|_3^2(t)+\int^t_0[\|W\|^2_3(0,\cdot,\tau)+(\|\r_t\|_1^2+\|\nabla\r\|^2_1+\|u\|^2_3)(\tau)]d\tau\\&\\&\leq
C\|W_0\|^2_3+C\d_0\int^t_0\|\nabla\r\|^2_2(\tau)d\tau.\ea \ee

From (\ref{3b}) it is easy to know that
$$
\|W_t\|_2^2(0,\cdot,t)\leq C\|W\|^2_3(0,\cdot,t),\ \
\|W_t\|_2^2(t)\leq C\|W\|_3^2(t),
$$
so (\ref{2417}) yields that
 \be\label{2418}\bal&\ \  \
\|W\|_3^2(t)+\|W_t\|_2^2(t)\\&\\&\ \ \ +\int^t_0[(\|W\|^2_3+\|W_t\|_2^2)(0,\cdot,\tau)+(\|\r_t\|_1^2+\|\nabla\r\|^2_1+\|u\|^2_3)(\tau)]d\tau\\
&\\&\leq C\|W_0\|^2_3+C\d_0\int^t_0\|\nabla\r\|^2_2(\tau)d\tau.
\ea\ee

By similar calculation to (\ref{2208}),
$\p^2(\ref{3b})_1\p^2\r_t+s\p^2(\ref{3b})_2\p^2\r_x+\p^2(\ref{3b})_3\p^2\r_y$
yields that
\be\label{2419}\bal&\int^t_0(\|\p^2\r_t\|^2+\|\p^2\r_y\|^2+\|\p^2\r_x\|^2)(\tau)d\tau\leq\\&\\&
C(\|W_0\|^2_3+(\|\nabla\r\|_2^2+\|u\|_2^2)(t)+
\int^t_0[(\|\r_t\|_2^2+\|W\|_2^2)(0,\cdot,\tau)+\|u\|_3^2(\tau)]d\tau)\\&\\&+C\d_0\int_0^t(\|u\|_3^2
+\|\nabla\r\|_2^2)(\tau)d\tau. \ea\ee

Choose $\lambda_9$ suitably small such that
$(\ref{2418})+\lambda_9(\ref{2419})$ yields that
 \be\label{2420}\bal&
\|W\|^2_3(t)+\|W_t\|_2^2(t)
+\int^t_0[(\|W\|^2_3+\|W_t\|_2^2)(0,\cdot,\tau)\\&\\&+(\|\r_t\|_2^2+\|\nabla\r\|_2^2+\|u\|^2_3)
(\tau)]d\tau\leq C\|W_0\|^2_3. \ea\ee

From (\ref{3b}) we know that $\|u_t\|_2^2(t)\leq
C(\|\nabla\r\|_2^2+\|u\|^2_3)(t)$, thus (\ref{2420}) yields that
 \be\label{2421}\bal&
\|W\|^2_3(t)+\|W_t\|_2^2(t)
+\int^t_0[(\|W\|^2_3+\|W_t\|_2^2)(0,\cdot,\tau)\\&\\&+(\|W_t\|_2^2+\|\nabla\r\|_2^2+\|u\|^2_3)
(\tau)]d\tau \leq C\|W_0\|^2_3.
 \ea\ee
\subsection{Estimate E}
By direct and a little tedious calculation, we get the following
estimates on the nonlinear terms.

\bl\label{25a} Assume (\ref{a}) holds, then
$$\bal&\ \
 |\int^t_0\int^{\infty}_{-\infty}\int^{\infty}_0\p_y^4H\cdot
\p_y^4Wdxdyd\tau|\\&\\&\leq
C\d_0(\|\p_y^4\r_0\|^2+\|\p_y^4\r\|^2(t)\\&\\&\ \
+\int^t_0[\|\p^4W\|^2(0,\cdot,\tau)+(\|\nabla\r\|^2_3+\|u\|^2_4)(\tau)]d\tau),
\ea
$$

$$\bal&\ \
 |\int^t_0\int^{\infty}_{-\infty}\int^{\infty}_0\p_y^3\p_xH\cdot
\p_y^3\p_xWdxdyd\tau|\\&\\&\leq
C\d_0(\|\p_y^3\p_x\r_0\|^2+\|\p_y^3\p_x\r\|^2(t)\\&\\&\ \
+\int^t_0[\|\p^4W\|^2(0,\cdot,\tau)+(\|\nabla\r\|^2_3+\|u\|^2_4)(\tau)]d\tau),
\ea
$$

$$\bal&\ \
 |\int^t_0\int^{\infty}_{-\infty}\int^{\infty}_0\p_y^2\p_x^2H\cdot
\p_y^2\p_x^2Wdxdyd\tau|\\&\\&\leq
C\d_0(\|\p_y^2\p_x^2\r_0\|^2+\|\p_y^2\p_x^2\r\|^2(t)\\&\\&\ \
+\int^t_0[\|\p^4W\|^2(0,\cdot,\tau)+(\|\nabla\r\|^2_3+\|u\|^2_4)(\tau)]d\tau),
\ea
$$

$$\bal&\ \
 |\int^t_0\int^{\infty}_{-\infty}\int^{\infty}_0\p_y\p_x^3H\cdot
\p_y\p_x^3Wdxdyd\tau|\\&\\&\leq
C\d_0(\|\p_y\p_x^3\r_0\|^2+\|\p_y\p_x^3\r\|^2(t)\\&\\&\ \
+\int^t_0[\|\p^4W\|^2(0,\cdot,\tau)+(\|\nabla\r\|^2_3+\|u\|^2_4)(\tau)]d\tau),
\ea
$$

$$\bal&\ \
 |\int^t_0\int^{\infty}_{-\infty}\int^{\infty}_0\p_x^4H\cdot
\p_x^4Wdxdyd\tau|\\&\\&\leq
C\d_0(\|\p_x^4\r_0\|^2+\|\p_x^4\r\|^2(t)\\&\\&\ \
+\int^t_0[\|\p^4W\|^2(0,\cdot,\tau)+(\|\nabla\r\|^2_3+\|u\|^2_4)(\tau)]d\tau).
\ea
$$
\el

As for the estimates on the boundary terms, we have the following
results.

\bl\label{25b} Assume (\ref{a}) holds, then
\be\label{2501}\|\p_y^3\p_xv_1\|^2(0,\cdot,t)\leq
C(\|W\|_3^2+\|\p_y^4u_2\|^2)(0,\cdot,t)+C\d_0\|W\|^2_4(0,\cdot,t),\ee
\be\label{2502}\bal \|\p_y^2\p_x^2v_1\|^2(0,\cdot,t)&\leq
C(\|W\|_3^2+\|\p_y^3\p_xu_2\|^2+\|\p_y^4(\r-u_1)\|^2)(0,\cdot,t)\\&\\&\
\ \ +C\d_0\|W\|^2_4(0,\cdot,t),\ea\ee \be\label{2503}\bal
\|\p_y\p_x^3v_1\|^2(0,\cdot,t)&\leq
C(\|W\|_3^2+\|\p_y^2\p_x^2u_2\|^2+\|\p_y^4u_2\|^2\\&\\
&\ \ \
+\|\p_y^3\p_x(\r-u_1)\|^2)(0,\cdot,t)+C\d_0\|W\|^2_4(0,\cdot,t),\ea\ee
\be\label{2504}\bal \|\p_x^4v_1\|^2(0,\cdot,t)&\leq
C(\|W\|_3^2+\|\p_y^3\p_xu_2\|^2+\|\p_y\p_x^3u_2\|^2+\|\p_y^2\p_x^2(\r-u_1)\|^2\\&\\
&\ \ \ +\|\p_y^4(\r-u_1)\|^2)(0,\cdot,t)
+C\d_0\|W\|^2_4(0,\cdot,t).\ea\ee
 \el

The proof of lemma \ref{25b} is similar to that of lemma \ref{23b},
so we omit here.

Multiplying $\p_y^4(\ref{3c})$ by $\p_y^4V$ and integrating it over
$\O_t$, we have that

\be\label{2505}\bal&\ \
\|\p_y^4W\|^2(t)+\int^t_0[\|\p_y^4W\|^2(0,\cdot,\tau)+\|\p^4_yu\|^2(\tau)]d\tau\\&\\&\leq
C(\|\p_y^4W_0\|^2+|\int^t_0\int^{\infty}_{-\infty}\int^{\infty}_0\p_y^4H\cdot
\p_y^4Wdxdyd\tau|). \ea\ee

Similarly, we have that
 \be\label{2506}\bal&\ \ \
\|\p_y^3\p_xW\|^2(t)\\&\\&\ \ \ +
\int^t_0[(\|\p_y^3\p_x(\r-u_1)\|^2+\|\p_y^3\p_xu_2\|^2)(0,\cdot,\tau)+\|\p_y^3\p_xu\|^2(\tau)]
d\tau\\&\\&\leq
C(\|\p_y^3\p_xW_0\|^2+\int^t_0\|\p_y^3\p_xv_1\|^2(0,\cdot,\tau)d\tau\\&\\&\
\ \  +|\int^t_0\int^{\infty}_{-\infty}\int^{\infty}_0
\p_y^3\p_xH\cdot \p_y^3\p_xWdxdyd\tau|), \ea\ee

\be\label{2507}\bal&\ \ \ \|\p_y^2\p_x^2W\|^2(t)\\&\\&\ \ \ +
\int^t_0[(\|\p_y^2\p_x^2(\r-u_1)\|^2+\|\p_y^2\p_x^2u_2\|^2)(0,\cdot,\tau)+\|\p_y^2\p_x^2u\|^2(\tau)]
d\tau\\&\\&\leq
C(\|\p_y^2\p_x^2W_0\|^2+\int^t_0\|\p_y^2\p_x^2v_1\|^2(0,\cdot,\tau)d\tau\\&\\&\
\ \  +|\int^t_0\int^{\infty}_{-\infty}\int^{\infty}_0
\p_y^2\p_x^2H\cdot \p_y^2\p_x^2Wdxdyd\tau|), \ea\ee

\be\label{2508}\bal&\ \ \ \|\p_y\p_x^3W\|^2(t)\\&\\&\ \ \ +
\int^t_0[(\|\p_y\p_x^3(\r-u_1)\|^2+\|\p_y\p_x^3u_2\|^2)(0,\cdot,\tau)+\|\p_y\p_x^3u\|^2(\tau)]
d\tau\\&\\&\leq
C(\|\p_y\p_x^3W_0\|^2+\int^t_0\|\p_y\p_x^3v_1\|^2(0,\cdot,\tau)d\tau\\&\\&\
\ \  +|\int^t_0\int^{\infty}_{-\infty}\int^{\infty}_0
\p_y\p_x^3H\cdot \p_y\p_x^3Wdxdyd\tau|), \ea\ee

and \be\label{2509}\bal&\ \ \ \|\p_x^4W\|^2(t) +
\int^t_0[(\|\p_x^4(\r-u_1)\|^2+\|\p_x^4u_2\|^2)(0,\cdot,\tau)+\|\p_x^4u\|^2(\tau)]
d\tau\\&\\&\leq
C(\|\p_x^4W_0\|^2+\int^t_0\|\p_x^4v_1\|^2(0,\cdot,\tau)d\tau+|\int^t_0\int^{\infty}_{-\infty}\int^{\infty}_0
\p_x^4H\cdot \p_x^4Wdxdyd\tau|). \ea\ee Choose
$\lambda_{10},~\lambda_{11},~\lambda_{12},~\lambda_{13}$ suitably
small such that $(\ref{2505})+\lambda_{10}(\ref{2506})$, combined
with (\ref{2421}) and (\ref{2501}), yields that
\be\label{2510}\bal&\ \
\|\p_y^4W\|^2(t)+\|\p_y^3\p_xW\|^2(t)+\int^t_0[(\|\p_y^4W\|^2+\|\p_y^3\p_x(\r-u_1)\|^2\\&\\&\
\  +\|\p_y^3\p_xu_2\|^2)(0,\cdot,\tau)
+(\|\p_y^4u\|^2+\|\p_y^3\p_xu\|^2)(\tau)] d\tau\\&\\&\leq
C(\|W_0\|_4^2+\d_0\int^t_0\|W\|_4^2(0,\cdot,\tau)d\tau\\&\\&\ \ +
|\int^t_0\int^{\infty}_{-\infty}\int^{\infty}_0\p_y^4H\cdot
\p_y^4Wdxdyd\tau|+|\int^t_0\int^{\infty}_{-\infty}\int^{\infty}_0
\p_y^3\p_xH\cdot \p_y^3\p_xWdxdyd\tau|), \ea\ee

 $(\ref{2510})+\lambda_{11}(\ref{2507})$, combined with (\ref{2421})
and (\ref{2502}), yields that
 \be\label{2511}\bal&\ \ \
(\|\p_y^4W\|^2+\|\p_y^3\p_xW\|^2+\|\p_y^2\p_x^2W\|^2)(t)\\&\\&\ \ \
+ \int^t_0[(\|\p_y^4W\|^2+\|\p_y^3\p_x(\r-u_1)\|^2+
\|\p_y^3\p_xu_2\|^2+ \|\p_y^2\p_x^2(\r-u_1)\|^2\\&\\&\ \ \
+\|\p_y^2\p_x^2u_2\|^2)(0,\cdot,\tau)
+(\|\p_y^4u\|^2+\|\p_y^3\p_xu\|^2+\|\p_y^2\p_x^2u\|^2)(\tau)]
d\tau\\&\\&\leq
C(\|W_0\|_4^2+\d_0\int^t_0\|W\|_4^2(0,\cdot,\tau)d\tau\\&\\&\ \ +
|\int^t_0\int^{\infty}_{-\infty}\int^{\infty}_0\p_y^4H\cdot
\p_y^4Wdxdyd\tau|+|\int^t_0\int^{\infty}_{-\infty}\int^{\infty}_0
\p_y^3\p_xH\cdot \p_y^3\p_xWdxdyd\tau|\\&\\&\ \ \ +
|\int^t_0\int^{\infty}_{-\infty}\int^{\infty}_0 \p_y^2\p_x^2H\cdot
\p_y^2\p_x^2Wdxdyd\tau|), \ea\ee

 $(\ref{2511})+\lambda_{12}(\ref{2508})$, combined with
(\ref{2421}) and (\ref{2503}), yields that

$$\bal&\ \ \
(\|\p_y^4W\|^2+\|\p_y^3\p_xW\|^2+\|\p_y^2\p_x^2W\|^2+\|\p_y\p_x^3W\|^2)(t)\\&\\&\
 \ \ +\int^t_0[(\|\p_y^4W\|^2+\|\p_y^3\p_x(\r-u_1)\|^2+
\|\p_y^3\p_xu_2\|^2 +\|\p_y^2\p_x^2(\r-u_1)\|^2\\&\\&\ \ \
+\|\p_y^2\p_x^2u_2\|^2+\|\p_y\p_x^3(\r-u_1)\|^2+\|\p_y\p_x^3u_2\|^2)(0,\cdot,\tau)\\&\\&\
\ \
+(\|\p_y^4u\|^2+\|\p_y^3\p_xu\|^2+\|\p_y^2\p_x^2u\|^2+\|\p_y\p_x^3u\|^2)(\tau)]
d\tau\ea$$

\be\label{2512}\bal&\leq
C(\|W_0\|_4^2+\d_0\int^t_0\|W\|_4^2(0,\cdot,\tau)d\tau\\&\\&\ \ \ +
|\int^t_0\int^{\infty}_{-\infty}\int^{\infty}_0\p_y^4H\cdot
\p_y^4Wdxdyd\tau|+|\int^t_0\int^{\infty}_{-\infty}\int^{\infty}_0
\p_y^3\p_xH\cdot \p_y^3\p_xWdxdyd\tau|\\&\\&\ \ \ +
|\int^t_0\int^{\infty}_{-\infty}\int^{\infty}_0 \p_y^2\p_x^2H\cdot
\p_y^2\p_x^2Wdxdyd\tau|\\&\\&\ \ \
+|\int^t_0\int^{\infty}_{-\infty}\int^{\infty}_0\p_y\p_x^3H\cdot
\p_y\p_x^3Wdxdyd\tau|), \ea\ee

and $(\ref{2512})+\lambda_{13}(\ref{2509})$, combined with
(\ref{2421}) and (\ref{2504}), yields that
 \be\label{2513}\bal&\ \ \
\|\p^4W\|^2(t)+\int^t_0[(\|\p_y^4W\|^2+\|\p_y^3\p_x(\r-u_1)\|^2+
\|\p_y^3\p_xu_2\|^2\\&\\&\ \ \
+\|\p_y^2\p_x^2(\r-u_1)\|^2+\|\p_y^2\p_x^2u_2\|^2 +
\|\p_y\p_x^3(\r-u_1)\|^2+\|\p_y\p_x^3u_2\|^2\\&\\&\ \ \
 +\|\p_x^4(\r-u_1)\|^2 +\|\p_x^4u_2\|^2)(0,\cdot,\tau) \ \
+\|\p^4u\|^2(\tau)] d\tau\\&\\&\leq
C(\|W_0\|_4^2+\d_0\int^t_0\|W\|_4^2(0,\cdot,\tau)d\tau+|\int^t_0\int^{\infty}_{-\infty}\int^{\infty}_0\p_y^4H\cdot
\p_y^4Wdxdyd\tau|\\&\\&\ \ \
+|\int^t_0\int^{\infty}_{-\infty}\int^{\infty}_0 \p_y^3\p_xH\cdot
\p_y^3\p_xWdxdyd\tau|\\&\\&\ \ \ +
|\int^t_0\int^{\infty}_{-\infty}\int^{\infty}_0 \p_y^2\p_x^2H\cdot
\p_y^2\p_x^2Wdxdyd\tau|\\&\\&\ \ \
+|\int^t_0\int^{\infty}_{-\infty}\int^{\infty}_0\p_y\p_x^3H\cdot
\p_y\p_x^3Wdxdyd\tau|\\&\\&\ \ \
+|\int^t_0\int^{\infty}_{-\infty}\int^{\infty}_0\p_x^4H\cdot
\p_x^4Wdxdyd\tau|). \ea\ee

Combined with lemma \ref{25a}, lemma \ref{25b} and (\ref{2421}),
(\ref{2513}) yields that \be\label{2515}\bal&\ \ \
\|W\|_4^2(t)+\int^t_0[\|W\|^2_4(0,\cdot,\tau)+(\|\r_t\|_2^2+\|\nabla\r\|^2_2+\|u\|^2_4)(\tau)]d\tau\\&\\&\leq
C\|W_0\|^2_4+C\d_0\int^t_0\|\nabla\r\|^2_3(\tau)d\tau.\ea \ee

From (\ref{3b}) it is easy to know that
$$
\|W_t\|_3^2(0,\cdot,t)\leq C\|W\|^2_4(0,\cdot,t),\ \
\|W_t\|_3^2(t)\leq C\|W\|_4^2(t),
$$
so (\ref{2515}) yields that
\be\label{2516}\bal&\ \  \
\|W\|_4^2(t)+\|W_t\|_3^2(t)\\&\\&\ \ \ +
\int^t_0[(\|W\|^2_4+\|W_t\|_3^2)(0,\cdot,\tau)+(\|\r_t\|_2^2+\|\nabla\r\|^2_2
+\|u\|^2_4)(\tau)]d\tau\\
&\\&\leq C\|W_0\|^2_4+C\d_0\int^t_0\|\nabla\r\|^2_3(\tau)d\tau.
\ea\ee

By similar calculation to (\ref{2208}),
$\p^3(\ref{3b})_1\p^3\r_t+s\p^3(\ref{3b})_2\p^3\r_x+\p^3(\ref{3b})_3\p^3\r_y$
yields that
\be\label{2517}\bal&\int^t_0(\|\p^3\r_t\|^2+\|\p^3\r_y\|^2+\|\p^3\r_x\|^2)(\tau)d\tau\leq\\&\\&
C(\|W_0\|^2_4+(\|\nabla\r\|_3^2+\|u\|_3^2)(t)+
\int^t_0[(\|\r_t\|_3^2+\|W\|_3^2)(0,\cdot,\tau)+\|u\|_4^2(\tau)]d\tau)\\&\\&+C\d_0\int_0^t(\|u\|_4^2
+\|\nabla\r\|_3^2)(\tau)d\tau. \ea\ee

Choose $\lambda_{14}$ suitably small such that
$(\ref{2516})+\lambda_{14}(\ref{2517})$ yields that
 \be\label{2518}\bal&
\|W\|^2_4(t)+\|W_t\|_3^2(t)
+\int^t_0[(\|W\|^2_4+\|W_t\|_3^2)(0,\cdot,\tau)\\&\\&+(\|\r_t\|_3^2+\|\nabla\r\|_3^2+\|u\|^2_4)
(\tau)]d\tau\leq C\|W_0\|^2_4. \ea\ee

From (\ref{3b}) we know that $\|u_t\|_3^2(t)\leq
C(\|\nabla\r\|_3^2+\|u\|^2_4)(t)$, thus (\ref{2518}) yields that
 \be\label{2519}\bal&
\|W\|^2_4(t)+\|W_t\|_3^2(t)
+\int^t_0[(\|W\|^2_4+\|W_t\|_3^2)(0,\cdot,\tau)\\&\\&+(\|W_t\|_3^2+\|\nabla\r\|_3^2+\|u\|^2_4)
(\tau)]d\tau\leq C\|W_0\|^2_4.
 \ea\ee

\subsection{Estimates on higher order derivatives}
By the similar arguments we can get the following estimates, for any
positive integer $l\geq4$ as long as $\d_0$ is sufficiently small,
 \be\label{2601}\bal&\ \ \
\|W\|^2_l(t)+\|W_t\|^2_{l-1}(t)\\&\\&\ \ \
+\int^t_0[(\|W\|^2_l+\|W_t\|^2_{l-1})
(0,\cdot,\tau)+(\|W_t\|^2_{l-1}+\|\nabla\r\|^2_{l-1}+\|u\|^2_l)(\tau)]d\tau\\&\\&\leq
C\|W_0\|^2_l.\ea\ee
\section{Theorems of existence}
\subsection{Local existence}
We are first going to obtain the local existence of solution to the
initial-boundary problem (\ref{3b}) by making use of iterative
scheme. Consider the following linear system,

\be\label{4a}
 \left\{ \bal&
\r^{m+1}_t-s\r^{m+1}_x+ru^{m+1}_{1x}+ru^{m+1}_{2y}=-r({\nabla\r^{m+1}\cdot
u^m+\r^m\rm{div}}u^{m+1}) ,
\\&\\&u^{m+1}_{1t}-su^{m+1}_{1x}+r\r^{m+1}_{x}+ku^{m+1}_1=-ru^{m}\cdot\nabla
 u^{m+1}_1+{1\over r}B^m\r^{m+1}_x,
\\&\\&u^{m+1}_{2t}-su^{m+1}_{2x}+r\r^{m+1}_{y}+ku^{m+1}_2=-ru^{m}\cdot\nabla u^{m+1}_2+
{1\over r}B^m\r^{m+1}_y,
\\&\\&(\r^{m+1},u^{m+1}_1,u^{m+1}_2)(x,y,t)|_{t=0}=(\r^{m+1}_0,u^{m+1}_{10},u^{m+1}_{20})(x,y),
\\&\\&(\r^{m+1}+u^{m+1}_1)|_{x=0}=0,
 \ea\right.
\ee where $B^m=r^2-{{P^{\prime} (1+\r^{m})}\over{1+\r^{m}}}$,
$\r^{m+1}_0,u^{m+1}_{10},u^{m+1}_{20}$ are functions of class
$C^{\infty}$ and
$\sum\limits_m\|\r^{m+1}_0-\r^{m}_0\|,\sum\limits_m\|u^{m+1}_{10}-u^{m}_{10}\|,\sum\limits_m\|u^{m+1}_{20}-u^{m}_{20}\|$
converge with the respective limits $\r_0, u_{10},u_{20}.$

Denote
$W^m=(\r^m,u^m_{1},u^m_{2}),W^m_0=(\r^m_0,u^m_{10},u^m_{20})$. By
the similar process to the a priori estimates in section 2, we have
the following estimate,

\be\label{4b}\bal&\ \ \ \|W^{m+1}\|^2_l(t)+\|W^{m+1}_t\|^2_{l-1}(t)+
\int^t_0[(\|W^{m+1}\|^2_l+\|W^{m+1}_t\|^2_{l-1})
(0,\cdot,\tau)\\&\\&+(\|W^{m+1}_t\|^2_{l-1}+\|\nabla\r^{m+1}\|^2_{l-1}+\|u^{m+1}\|^2_l)(\tau)]d\tau\\&\\&\leq
C(\|W^m_0\|_l)\|W^{m+1}_0\|^2_l+C(\|W^m\|_l(t))\|\nabla\r^{m+1}\|^2_{l-1}(t)\\&\\&+
\int^t_0C(\|W^m\|_l(\tau),\|W^{m+1}\|_l(\tau))[\|W^{m+1}\|^2_l
(0,\cdot,\tau)\\&\\&+(\|\nabla\r^{m}\|^2_{l-1}+\|u^{m}\|^2_l+\|\nabla\r^{m+1}\|^2_{l-1}+\|u^{m+1}\|^2_l)(\tau)]d\tau.
\ea\ee

From (\ref{4b}) we get the following lemma for the system
(\ref{4a}).

\bl\label{411} Let $l$ be an integer, $l\geq4.$ Assume that
$\r_0,u_{10},u_{20}\in H^l(\R^+\times\R),$ and
$\|\r_0\|_{l},\|u_{10}\|_{l},\|u_{20}\|_{l}$ are sufficiently small.
Then there exists a time $T_1$ and a number $R_1$, such that for all
$m\geq0,$ we have
$$
\sup\limits_{0\leq t\leq T_1}\|W^m\|_l(t)\leq R_1, \ \
\sup\limits_{0\leq t\leq T_1}\|\partial_tW^m\|_{l-1}(t)\leq R_1,
$$
where the numbers $R_1, T_1$ depend both on the system (\ref{4a})
and on the initial data
$\|\r_0\|_{l},\|u_{10}\|_{l},\|u_{20}\|_{l}$. \el

Now we are going to show the convergence of the iterative scheme in
$L^2(\R^+\times\R)$ on a smaller time interval $T^{\ast}$, then we
conclude the convergence in $H^r(\R^+\times\R)$ for all $0\leq r<l$
by interpolation.

First we define the difference $\bar{W}^m\triangleq W^{m+1}-W^m$ and
other denotations can be similarly defined. We form the difference
of two successive equations of the scheme,

\be\label{4c}
 \left\{ \bal&
\bar{\r}^{m}_t-s\bar{\r}^{m}_x+r\bar{u}^{m}_{1x}+r\bar{u}^{m}_{2y}=\bar{h}^m_1
,
\\&\\&\bar{u}^{m}_{1t}-s\bar{u}^{m}_{1x}+r\bar{\r}^{m}_{x}+k\bar{u}^{m}_1=\bar{h}^m_2,
\\&\\&\bar{u}^{m}_{2t}-s\bar{u}^{m}_{2x}+r\bar{\r}^{m}_{y}+k\bar{u}^{m}_2=\bar{h}^m_3,
\\&\\&(\bar{\r}^{m},\bar{u}^{m}_1,\bar{u}^{m}_2)(x,y,t)|_{t=0}=(\bar{\r}^{m}_0,\bar{u}^{m}_{10},\bar{u}^{m}_{20})(x,y),
\\&\\&(\bar{\r}^{m}+\bar{u}^{m}_1)|_{x=0}=0,
 \ea\right.
\ee where

$$\bal\bar{h}^m_1&=-r(\nabla\r^{m+1}\cdot
u^m+\r^m{\rm{div}}u^{m+1}-\nabla\r^{m}\cdot
u^{m-1}-\r^{m-1}{\rm{div}}u^{m})\\&\\&=-r(\nabla\bar{\r}^{m}\cdot
u^m+{\r^m\rm{div}}\bar{u}^{m}+\nabla\r^{m}\cdot
\bar{u}^{m-1}+\bar{\r}^{m-1}{\rm{div}}u^{m}), \ea$$

$$\bal\bar{h}^m_2&=-ru^{m}\cdot\nabla
 u^{m+1}_1+{1\over r}B^m\r^{m+1}_x+ru^{m-1}\cdot\nabla
 u^{m}_1-{1\over r}B^{m-1}\r^{m}_x\\&\\&=-r(u^{m}\cdot\nabla
 \bar{u}^{m}_1+\bar{u}^{m-1}\cdot\nabla
 u^{m}_1)+{1\over r}B^m\bar{\r}^{m}_x+{1\over r}\bar{B}^{m-1}\r^{m}_x,
 \ea$$

$$\bal\bar{h}^m_3&=-ru^{m}\cdot\nabla
 u^{m+1}_2+{1\over r}B^m\r^{m+1}_y+ru^{m-1}\cdot\nabla
 u^{m}_2-{1\over r}B^{m-1}\r^{m}_y\\&\\&=-r(u^{m}\cdot\nabla
 \bar{u}^{m}_2+\bar{u}^{m-1}\cdot\nabla
 u^{m}_2)+{1\over r}B^m\bar{\r}^{m}_y+{1\over r}\bar{B}^{m-1}\r^{m}_y.
 \ea$$

 By the similar process to the a priori estimates in section 2, we
 get the following estimate for the system (\ref{4c}),

 \be\label{4d}\bal
&\|\bar{W}^m\|^2(t)\\&\\&\leq
C_1(R_1)\int^t_0\|\bar{W}^{m-1}\|^2(\tau)d\tau+C_2(R_1)\int^t_0\|\bar{W}^{m}\|^2(\tau)d\tau+C_3(R_1)\|\bar{W}^{m}_0\|^2.
 \ea\ee
Denote $$y_m\triangleq\sup\limits_{0\leq t\leq
T^{\ast}}\|\bar{W}^m\|^2(t),$$ then we have from  (\ref{4d}),
$$y_m\leq C_2(R_1)T^{\ast}y_m+C_1(R_1)T^{\ast}y_{m-1}+\beta_m,$$
where $\beta_m=C_3(R_1)\|\bar{W}^{m}_0\|^2.$ We choose $T^{\ast}$ to
be such that $$(C_1(R_1)+C_2(R_1))T^{\ast}\leq {1\over2}.$$ It
yields that

\be\label{4e}\sum\limits_my_m\leq2\sum\limits_m\beta_m.\ee
 By using
lemma 3.6.5 in \cite{S}(see page 98), we know that
$\{\beta_m\}_{m\geq0}$ has a finite sum. From (\ref{4e}) we deduce
that $\{y_m\}_{m\geq0}$ equally has a finite sum, that is to say
$W^m$ converges at least in
$L^{\infty}([0,T^{\ast}];L^2(\R^+\times\R))$. We denote the limit as
$W=(\r,u_1,u_2)$, then $W\in
L^{\infty}([0,T^{\ast}];L^2(\R^+\times\R))$. By an interpolation
formula between $H^0=L^2$ and $H^l$, we have for all $0\leq r<l$,
$$
\|W^m-W\|_r\leq|W^m-W|_2^{1-{r\over l}}\|W^m-W\|_l^{r\over l}.
$$
So the sequence $\{W^m\}_{m\geq0}$ tend to $W$ in
$L^{\infty}([0,T^{\ast}];H^r(\R^+\times\R))$ for all $r<l$. Since
$l\geq4,$ we have the result that $W$ is a regular solution of the
initial-boundary value problem (\ref{3b}). So we obtain the
following theorem of local existence.

\bt\label{412} Let $l$ be an integer, $l\geq4.$ Assume that
$\r_0,u_{10},u_{20}\in H^l(\R^+\times\R),$ and
$\|\r_0\|_{l},\|u_{10}\|_{l},\|u_{20}\|_{l}$ are sufficiently small.
Then there exists a time $T>0$ such that the problem (\ref{3b}) has
a unique classical solution $$(\r,u_1,u_2)\in
C^1([0,T]\times\R^+\times\R).$$ In addition, $(\r,u_1,u_2)\in
C^1([0,T];H^{l-1}(\R^+\times\R))\cap C^0([0,T];H^l(\R^+\times\R)).$
\et

{\bf Remark.} As mentioned in the section of the introduction, for
the initial-boundary value problem to the isentropic Euler equations
with damping, we obtain the local existence of the classical
solution only in the case of the small initial data due to some
essential or technical difficulties, while for the Cauchy problem of
symmetric hyperbolic systems, the local existence of classical
solutions  can be proved by using the fixed point mapping theorem or
the iteration method without the assumption that the initial data
are small (see \cite{S}).
 \subsection{Global existence}In order to obtain
the global existence of classical solution to the system (\ref{3b}),
we only need to prove the a priori estimate. Based on the preceding
estimates in section 2, (\ref{2601}) yields the a priori assumption
(\ref{a}) for any time $T$. Therefore we have the following theorem
of global existence.

\bt\label{result} Assume that $\r_0,u_{10},u_{20}\in
H^l(\R^+\times\R),$ $l\geq4$ is a positive integer, and
$\|\r_0\|_{l},\|u_{10}\|_{l},\|u_{20}\|_{l}$ are sufficiently small.
Then there exists a unique, global, classical solution
$(\r,u_1,u_2)$ to the initial-boundary value problem (\ref{3b})
which satisfies (\ref{2601}) and $$ (\r,u_1,u_2)\in
C^1([0,\infty);H^{l-1}(\R^+\times\R))\cap
C^0([0,\infty);H^{l}(\R^+\times\R)).
$$\et

{\bf Remark.} 1. In this paper, although we study the IBVP for 2-D
Euler equations with damping, in fact the corresponding results
still hold in the case of $n$-D $(n\geq3)$.

2. In this paper, we assume that the boundary function in (\ref{3a})
is constant, and it results in the homogeneous boundary condition in
(\ref{3b}), so the estimates of the solution can be controlled only
by the initial data, otherwise they should be controlled by both the
initial and the  boundary functions.


\begin{thebibliography}{}
\bibitem{Kreiss}H. O. Kreiss and J. Lorenz, Initial-boundary value problems
 and the Navier-Stokes equations, Classics
 in Applied Mathematics, 47. Society for Industrial and Applied Mathematics (SIAM),
 Philadelphia, PA, 2004.
\bibitem{HMS} F. M. Huang, A. Matsumura and X. D. Shi, Viscous shock
wave and boundary layer solution to an inflow problem for
compressible viscous gas. Comm. Math. Phys.
 239 (2003), no. 1-2, 261--285.
\bibitem{KK1}Y. Kagei and T. Kobayashi, On large-time behavior
of solutions to the compressible Navier-stokes equations in the half
space in $\R^3$, Arch. Rational Mech. Anal. 165 (2002), 89¨C159.
 \bibitem{KK}Y. Kagei and S. Kawashima, Stability of planar stationary solutions
to the compressible Navier-Stokes equation on the half space,
Commun. Math. Phys. 266 (2006), 401¨C430.
 \bibitem{MM} A. Matsumura and M. Mei, Convergence to travelling fronts of solutions of the $p$-system with
 viscosity in the presence of a boundary. Arch. Ration. Mech. Anal. 146 (1999), no. 1, 1--22.
 \bibitem{MN1}A. Matsumura and K. J. Nishihara, Large-time behaviors of solutions
 to an inflow problem in the half space for a one-dimensional system of compressible
 viscous gas. Comm. Math. Phys. 222 (2001), no. 3, 449--474.
 \bibitem{MN2}A. Matsumura and K. J. Nishihara, Global asymptotics toward the rarefaction wave
 for solutions of viscous $p$-system with boundary effect. Quart. Appl. Math. 58 (2000), no. 1, 69--83.
 \bibitem{NY}K. J. Nishihara and T. Yang, Boundary effect on asymptotic
behavior of solutions to the p-system with linear damping, J.
Differential Equations 156 (1999), 439-458.
 \bibitem{S} D. Serre, Systems of
conservation laws I, hyperbolicity entropy, shock waves, Cambridge
University Press (1999).
\bibitem{WY}W. K. Wang and T. Yang, The pointwise estimates of
solutions  for Euler equations with damping in multi-dimensions, J.
Differential Equations 173 (2001), 410-450.



\end{thebibliography}
\end{document}